\documentclass[11pt, a4paper]{article}
\usepackage{amsmath,amssymb,latexsym,graphics,epsfig}
\usepackage{hyperref}
\usepackage{color}
\usepackage{amsthm}

\input amssym.def
\newsymbol\rtimes 226F
\newfont{\nset}{msbm10}
\newcommand{\ns}[1]{\mbox{\nset #1}}

\def\R{\ns{R}}

\def\A{\mbox{\boldmath $A$}}

\def\D{\mbox{\boldmath $D$}}

\def\J{\mbox{\boldmath $J$}}

\def\L{\mbox{\boldmath $L$}}

\def\R{\ns{R}}

\def\diag{\mathop{\rm diag}\nolimits}

\def\>{\mathop{\rightarrow}\nolimits}

\def\f{\mbox{\boldmath $f$}}

\def\j{\mbox{\boldmath $j$}}

\def\vecdelta{\mbox{\boldmath $\delta$}}

\def\vecnu{\mbox{\boldmath $\nu$}}

\def\vec0{\mbox{\bf 0}}

\def\tr{\mathop{\rm tr}\nolimits}

\begin{document}

\title{Number of Walks and \\ Degree Powers in a Graph}
\author{M.A. Fiol, E. Garriga
\\ \\
{\small Universitat Polit\`ecnica de Catalunya, BarcelonaTech}\\
{\small Departament de Matem\`atica Aplicada IV} \\
{\small Barcelona, Catalonia} \\
{\small (e-mails: {\tt \{fiol,egarriga\}@ma4.upc.edu})}}

\date{}

\maketitle

\begin{abstract}
This note deals with the relationship between the total number of
$k$-walks in a graph, and the sum of the $k$-th powers of its vertex
degrees. In particular, it is shown that the the number of all
$k$-walks is upper bounded by the sum of the $k$-th powers of the
degrees.
\end{abstract}


Let $G =(V,E)$ be a connected graph on $n$ vertices,
$V=\{1,2,\ldots,n\}$, with adjacency matrix $\A$. For any integer
$k\ge 1$, let $a_{ij}^{(k)}$ denote the $(i,j)$ entry of the power
matrix $\A^k$. Let $\D$ be the diagonal matrix with elements
$(\D)_{ii} = d_i$ (the degree of vertex $i$). Here we study the
relationship between the sum of all walks of length $k$ in $G$ and
the sum of the $k$-th powers of its degrees. As a main result, and
answering in the affirmative a conjecture of Marc Noy [8], we will
show that
\begin{equation}
\label{basic-ineq} \sum_{i,j} a_{ij}^{(k)}\le \sum_{i}d_i^k,
\end{equation}
with equality if and only if $G$ is regular or $k\le 2$. In the case
$k=3$ we also provide an exact value of the difference between the
above sums in (\ref{basic-ineq}). In other line of research, some
upper bounds for $\sum_{i}d_i^k$ have been given by several authors.
See, for instance, \cite{c06,d98,dc98,pps99} (for general graphs)
and \cite{bn04,cy00} (for graphs not containing a prescribed
subgraph).

Let us first begin with the small values of $k$. The case $k=0$ is
trivial since the number of walks of length $0$ equals the number of
vertices. Similarly, if $k=1$, the sum $\sum_{i,j} a_{ij}$ is just
the sum of the degrees $d_1+d_2 +\cdots + d_n$. If $k=2$, we can use
that $\A\j = \D\j$ (where $\j$ is the all-$1$ vector) and the
symmetry of the involved matrices to obtain:
$$
\sum_{i,j} a_{ij}^{(2)}= \langle\j,\A^2\j\rangle=\langle
\A\j,\A\j\rangle= \|\A\j\|^2= \|\D\j\|^2 = d_1^2 + d_2^2 +\cdots +
d_n^2.
$$

Assume now that $G$ is regular of degree, say, $d$. Then, $\j$ is
the positive eigenvector corresponding to the eigenvalue $d$ and we
get
$$
\sum_{i,j}a_{ij}^{(k)}=\langle \j, \A^k\j\rangle = \langle
\j,d^k\j\rangle = d^k\|\j\|^2 = nd^k.
$$
A similar reasoning shows that, for a general (non-regular) graph,
the inequality in (\ref{basic-ineq}) always holds for $k$ large
enough. Indeed, let $\vecnu =(\nu_1,\nu_2,\ldots,\nu_n)^{\top}$ be
the positive eigenvector of $G$, normalized in such a way that
$\min_{i\in V}\nu_i= 1$. Let $\lambda$ be its corresponding
(positive) eigenvalue, which is known to be smaller than the maximum
degree $\Delta$ of $G$ (see \cite{b93,c71}). Then,
$$
\sum_{i,j}a_{ij}^{(k)}\le \langle \vecnu, \A^k\vecnu \rangle =
\langle \vecnu, \lambda^k\vecnu \rangle = \|\vecnu\|^2\lambda^k
$$
which, for $k$ large enough, is smaller than the single $k$-power
$\Delta^k$.

To deal with the case $k=3$, it is more convenient to work with the
Laplacian matrix of $G$; that is, $\L :=\D-\A$. Then, recall that,
given any real function defined on $V$, $f:V\rightarrow \R$, and
with $\f$ being the (column) vector with components the values of
$f$ on $V$, we have
$$
\langle \f, \L\f\rangle = \sum_{i\sim j}(f(i) - f(j))^2,
$$
where the sum is extended over all edges of $G$ (see, for instance,
\cite{b93}). We are interested in the case when the above function
is just the degree of the corresponding vertex: $f(i)= d_i$. In this
case, we denote by $\vecdelta$ its corresponding vector. Then, the
difference between the two sums in (\ref{basic-ineq}) is just
$$
\sum_{i}d_i^3-\sum_{i,j}a_{ij}^{(3)}=\langle \vecdelta,
\D\vecdelta\rangle-\langle \vecdelta, \A\vecdelta\rangle=\langle
\vecdelta, \L\vecdelta\rangle=\sum_{i\sim j} (d_i-d_j)^2\ge 0.
$$

To prove the inequality in the general case, note first that, for
any two positive numbers $a,b$ with, say, $a\ge b$,
$$
a^rb + ab^r = a^{r+1} + b^{r+1} - (a^r - b^r)(a - b) \le a^{r+1} +
b^{r+1}
$$
with equality if and only if $a = b$. (The same conclusion is
reached when we apply H\"{o}lder inequality \cite{hlp34} to the
vectors $(a,b)$ and $(b^r,a^r)$ with norms $L^{r+1}$ for the first
vector and dual norm $L^{(r+1)/r}$ for the second one.)

Also, notice that all walks of a given length, say $k\ge 1$, can be
obtained by considering, for any vertex $i$, all $(i,j)$-walks of
length $k-1$ ``extended" by each of the $d_i$ edges incident to $i$.
Thus, $\sum_{i,j} a_{ij}^{(k)} = \sum_{i,j} d_i a_{ij}^{(k-1)}$ or,
equivalently,
$$
\langle \j, \A^k\j \rangle = \langle \j, \A^{k-1}\D\j\rangle =
\langle \D\j, \A^{k-1}\j\rangle.
$$
Keeping all this in mind, we are now ready to prove
(\ref{basic-ineq}). Indeed, assuming that $k\ge 3$, we have:
\begin{eqnarray*}
\sum_{i,j} a_{ij}^{(k)} & = & \sum_{i,j} d_i a_{ij}^{(k-2)} d_j =
\sum_{i} a_{ii}^{(k-2)} d_i^2 + \sum_{i<j} 2a_{ij}^{(k-2)} d_i d_j \\
   & \le & \sum_{i} a_{ii}^{(k-2)} d_i^2 + \sum_{i<j} a_{ij}^{(k-2)}
(d_i^2+ d_j^2) \\
   & = & \sum_{i,j} a_{ij}^{(k-2)} d_j^2 \\
   & = & \sum_{i,j} d_i a_{ij}^{(k-3)} d_j^2 = \sum_{i} a_{ii}^{(k-3)}d_i^3
   + \sum_{i<j} a_{ij}^{(k-3)} (d_id_j^2+ d_i^2 d_j) \\
   & \le & \sum_{i} a_{ii}^{(k-3)}d_i^3
   + \sum_{i<j} a_{ij}^{(k-3)} (d_i^3+ d_j^3) \\
   & = & \sum_{i,j} a_{ij}^{(k-3)}d_j^3 \le \cdots \le \sum_{i,j} a_{ij}
   d_j^{k-1} = \sum_j d_j^k.
\end{eqnarray*}
Moreover, notice that all the above inequalities become equalities
if and only if $G$ is regular, as claimed.

As commented by one of the referees, the above prove suggests a
wider result: If $\A$ is a real symmetric $n\times n$ matrix,
$\D=\diag (d_1,d_2,\ldots,d_n)$, where $d_i=\sum_{j}|a_{ij}|$, and
$\J$ is the all-$1$ matrix, then
\begin{equation}
\label{basic-ineq-matrix}
\tr (\A^k\J) \le \tr (\D^k\J)=\sum_{i}
d_i^k.
\end{equation}
In other words, an extremal (not characteristic) property of the
positive diagonal matrices is exhibited.

\subsection*{Acknowledgments}
Research supported by the Ministerio de Educaci\'on y Ciencia,
Spain, and the European Regional Development Fund under project
MTM2005-08990-C02-01 and by the Catalan Research Council under
project 2005SGR00256.

The authors thank one of the referees for the comment about the use
of H\"older inequality and the result in (\ref{basic-ineq-matrix}).


\begin{thebibliography}{99}
\bibitem{b93}
N. Biggs, {\em Algebraic Graph Theory}, second ed., Cambridge
University Press, Cambridge, 1993.

\bibitem{bn04}
B. Bollob\'as and V. Nikiforov, Degree powers in graphs with
forbidden subgraphs, {\em  Electron. J. Combin.} 11 (2004), Research
paper R42.

\bibitem{c06}
Sebastian M. Cioab\v{a}, Sums of powers of the degrees of a graph,
{\em Discrete Math.} 306 (2006) 1959--1964.

\bibitem{c71}
D.M. Cvetkovi\'c, Graphs and their spectra, {\em Univ. Beograd.
Publ. Elektrotehn. Fak. Ser. Mat. Fiz.} No. 354--356 (1971), 1--50.

\bibitem{cy00}
Y. Caro and R. Yuster, A Tur\'an type problem concerning the powers
of the degrees of a graph, {\em Electron. J. Combin.} 7 (2000),
Research paper R47.

\bibitem{d98}
K. Ch. Das, Maximizing the sum of the squares of the degrees in a
graph, {\em Discrete Math.} 285 (2004), 57--66.

\bibitem{dc98}
D. de Caen, An upper bound on the sum of the squares of the degrees
in a graph, {\em Discrete Math.} 185 (1998), 245--248.

\bibitem{hlp34}
G.H. Hardy, J.E. Littlewood and G. P\'olya,  {\em Inequalities},
Cambridge Univ. Press, Cambridge, 1934.

\bibitem{n05}
Marc Noy, personal communication (2005).

\bibitem{pps99}
U.N. Peled, R. Petreschi and A. Sterbini, $(n,e)$-graphs with
maximum sum of squares of degrees, {\em J. Graph Theory} 31 (1999),
283--295.
\end{thebibliography}
\end{document}